\documentclass[times]{elsarticle}
\usepackage{graphicx,colordvi}
\begin{document}
\begin{frontmatter}

\title{Precision of the ENDGame: Mixed-precision arithmetic in the iterative solver of the Unified Model.}

\author[met,read]{C.~M.~Maynard}
\ead{c.m.maynard@reading.ac.uk}
\author[met]{D.~N.~Walters}

\address[met]{Met Office, FitzRoy Road, Exeter, EX1 3PB}
\address[read]{Department of Computer Science, Polly Vacher Building, University of Reading, Reading, UK, RG6 6AY}

\begin{abstract}
  The Met Office's weather and climate simulation code the Unified
  Model is used for both operational Numerical Weather Prediction and
  Climate modelling. The computational performance of the model
  running on parallel supercomputers is a key consideration. A Krylov sub-space solver is employed
  to solve the equations of the dynamical core of the model, known as
  ENDGame. These describe the
  evolution of the Earth's atmosphere. Typically, 64-bit precision is
  used throughout weather and climate applications. This work presents
  a mixed-precision implementation of the solver, the beneficial effect on
  run-time and the impact on solver convergence. The complex interplay of errors arising from
  accumulated round-off in floating-point arithmetic and other
  numerical effects is discussed. A careful analysis is required, however,
  the mixed-precision solver is now employed in the operational forecast to
  satisfy run-time constraints without compromising the accuracy of
  the solution. 
\end{abstract}


\end{frontmatter}

\section{Introduction}
\label{sec:intro}

Numerical simulations of the fluid dynamics of the Earth's atmosphere,
which are central to all weather and climate models, require
solutions to a non-linear system of partial differential
equations (PDEs). These PDEs are recast into a discrete system of
algebraic equations on a grid covering the Earth. For each grid cell,
the fluid dynamic properties of the atmosphere, the velocity,
temperature, pressure and density of the air, are then to be
determined. 

When deriving these algebraic equations, which arise from
the continuous form of the compressible Euler equations, there is a
choice of whether to integrate forward in time explicitly (a direct
calculation) or (semi-)implicitly. The latter involves a global matrix
inversion of some form. For Numerical Weather Prediction (NWP)
simulations, particularly those discretised onto a longitude-latitude (lon-lat) grid, 
the nature of the grid means that it is generally
not feasible to use an explicit scheme as it would severely restrict the
length of the time step employed in the simulation.
A typical operational global NWP simulation is required to simulate seven to fourteen 
days within one-to-two hours of wall-clock time and the number of small
time steps of the explicit scheme would either be prohibitively
expensive, or severely limit the physical complexity or spatial
resolution of an affordable simulation.

Semi-implicit schemes treat the short-timescale accoutic-gravity
modes implicitly and allow a stable integration of the flow around the
polar singularity and in regions with strong horizontal flows. 
To an extent, these schemes also allow more flexibility in the
formulation of the algebraic system and in particular the form of the
global matrix inversion. In line with incompressible and low Mach
number flows, a pressure correction equation is derived which takes
the form
\begin{equation}
A\cdot{\mathbf x} = {\mathbf b},
\label{eqn:LinearSys}
\end{equation}
where $A$ is a large, sparse (banded) matrix of order $n$, ${\mathbf x}$ is
the pressure correction / tendency and ${\mathbf b}$ contains the explicit forcing
terms. This can be solved using a combination of Krylov sub-space 
iterative solvers, pre-conditioners and parallel computers to satisfy
wall-clock time constraints such as those described above.

The current operational configuration of the Met Office Unified Model
(UM)~\cite{Brown:2012} solves the dynamical system described above
using the ``ENDGame'' dynamical core~\cite{QJ:QJ2235}. The operational 
deterministic global forecast is discretised onto an N1280 lon-lat grid with
$2560\times1920$ grid points, giving it a horizontal resolution of
approximately 10\,km at mid-latitudes. With 70 discrete vertical levels,
the number of degrees of freedom $n$ for the pressure correction
system is large at approximately 350 million; the full algebraic system
is six times larger. For such a large system, a parallel supercomputer
is required. However, even on state-of-the-art parallel computers
such as the Met Office's Cray XC40 machine, both algorithmic and code
optimisations are necessary to scale to a large degree of
parallelism\footnote{The operational global forecast model uses 
  around 500 compute nodes, each with 36 CPU cores.} and to execute the
program as quickly as possible.   

One potential optimisation is to consider the numerical precision
to which variables in the code are computed and stored, as 
reducing the precision can provide performance advantages.
In common with many scientific numerical applications, weather and
climate codes are memory bound. Using smaller data types reduces data movement between memory
and CPU and reduces remote communication between processors, such that this is faster and uses less
energy. Moreover, caches can be utilised more effectively. 
For weather and climate models, however, knowledge of accuracy and
uncertainty are important.
There are many processes that affect the evolution of the system that are either
not resolved at the cut-off scale imposed by the discretisation or are
not represented by the dynamical core. These processes are represented
by physical parametrizations, which often include many branches in the
code. For example, the presence or not of liquid water in a grid cell
will determine whether latent heat can be released through the
freezing of that water;
this in turn will affect the temperature of the atmosphere and so
influence the evolution of the dynamics.
An aversion to the possibility of compromising the accuracy of a
simulation through either accumulated round-off errors, or round-off
errors triggering different branches with parametrizations, has led to the
``safety first'' approach of using high precision arithmetic throughout, and the
assumption that this is both necessary and more accurate. 
In common with many scientific
applications, therefore, 
64-bit floating point arithmetic -- often referred to as double precision
-- is used by default and has become the accepted standard.

Recent work~\cite{Palmer_SP, DUBEN20142, Dueben2015, Vana2017}
and~\cite{Nakano_SP} has explored the potential benefits of moving
away from the use of 64-bit arithmetic throughout an atmospheric model.
These papers examine a variety of approaches to reducing precision,
from low precision emulators, Field Programmable Gate Arrays (FPGAs),
toy models, right through to running a full atmospheric model in
32-bit precision ({\em i.e.} single precision) to examine the effect of precision on
forecast accuracy. Whilst it is true that 64-bit precision may be
desirable, if not required, for some portions of the simulation, many
schemes within a model will have either physical or parametric
uncertainty, or will be designed to reach an approximate solution,
with accuracy far removed from the numerical precision of the
variables within the scheme. 
One such example is the numerical solution of
Equation~(\ref{eqn:LinearSys}). The solution for the
pressure correction $\mathbf x$ is obtained through an iterative
solver, in which each variable is typically defined to double
precision, but in which the calculation is halted when the normalised
error has reduced by only a small number of orders of magnitude.
Whilst the solver is only a small proportion of the total code base
when measured by the number of lines of code, at the resolutions and
node counts described for operational global NWP above, these routines
constitute approximately one quarter of the run time of the full simulation.
In this work, a mixed-precision solver is defined, in which the
pressure field itself is held in double precision throughout, but the
majority of calculations, and hence the majority of memory transfers,
are performed in single precision.
In section~\ref{sec:ENDGame} the ENDGame dynamical core,
the Helmholtz equation to be solved for the pressure correction and
the pre-conditioner used for this are described. In section~\ref{sec:MPSolver} the
mixed-precision solver is described, whilst in section~\ref{sec:res}
the effect this has on convergence to the solution as
well as the computational speed up obtained is discussed. Finally, conclusions are
drawn in section~\ref{sec:con}. 
 
\section{The ENDGame Dynamical core}
\label{sec:ENDGame}
Symbolically, the linear system of equations to be solved is of the form
\begin{equation}
   \left( \begin{array}{cccccc}
          I & 0 & 0 & 0 & 0 & G_u \\
          0 & I & 0 & 0 & 0 & G_v \\
          0 & 0 & I & B_i & 0 & G_w \\
          0 & 0 & B_{ii} & I & 0 & 0 \\
        D_x & D_y & D_z & 0 & I & 0 \\
          0 & 0 & 0 & E_\theta & I & E_\pi 
          \end{array}
   \right) \left( \begin{array}{c} u \\ v \\ w \\ \theta \\ \rho \\ \Pi \end{array} \right)
   =  \left( \begin{array}{c} R_u \\ R_v \\ R_w \\ R_\theta \\ R_\rho \\ R_\Pi \end{array} \right)
\end{equation}
where $I$ represents identity matrices, $G_{u,v,w}$ are discrete
pressure gradient terms, $D_{x,y,z}$ are components of the divergence
operator, $B_i$ and $B_{ii}$ are the coupling between the potential
temperature and the vertical component of velocity due to gravity and
$E_{\theta, \rho}$ arise from linearisation of the equation of
state. The solution of this system corresponds to the fluid velocity
components $(u,v,w)$ and the tendencies (change between time steps) of
the thermodynamic variables $\rho$ (density), $\Pi$ (Exner pressure)
and $\theta$ (potential temperature\footnote{The absolute temperature
  $T = \theta \times \Pi$.}). Note that within the ENDGame
formulation, due to non-linearities and the change to a terrain
following coordinate system, the right hand side terms involve
contributions from the previous iterates of the current time
step. Typically this linear system is solved four times per time step
and so the solution method needs to be efficient.

This matrix is of the general form
\begin{equation}
  \left( \begin{array}{cc}
           P & Q \\
           R & S 
   \end{array}
   \right) \left( \begin{array}{c} {\mathbf y} \\ {\mathbf x} \end{array} \right)
   =  \left( \begin{array}{c} {\mathbf s} \\ {\mathbf t} \end{array} \right)
\end{equation}
which can be solved using the Schur compliment as
\begin{equation}
   {\mathbf y} = P^{-1}( {\mathbf s}- Q{\mathbf x} ); \hspace{0.1in}
   A{\mathbf x}  = ( S - RP^{-1}Q) {\mathbf x} = {\mathbf t} -
   RP^{-1}{\mathbf s} = {\mathbf b}.
  \label{eqn:Helmholtz}
\end{equation}
The matrix $A$ has a
$7$-band structure similar to what would be expected from discretising
the Laplacian operator, which, would be its exact equivalent in the
Boussinesq incompressible limit. In the full system, however, it differs in some key
aspects. Firstly, the matrix is non-symmetric and it is not constant
coefficient; {\em i.e.}, the matrix rows are all different due to the
spherical nature of the Earth and variations in temperature and orography
(the height of the lower boundary). Secondly, as may be anticipated from the previous
sentence, the matrix is time varying and so needs to be recomputed at
every time step. This precludes the ability to pre-factorise the
matrix off-line. Finally, because the vertical length scales are much
shorter than the horizontal and gravitational effects are far larger
than the vertical accelerations, some care is needed in order to solve
this matrix. In particular, it is necessary to ensure that the large
scale hydrostatic balance is maintained during the iterations. It
should be noted that, since this equation arises from a non-linear
algebraic system, there are terms involving ${\mathbf x}$ which are lagged
and so appear as forcing in the right-hand-side ${\mathbf b}$
(see~\cite{QJ:QJ2235} for further details).

The system in equation~(\ref{eqn:Helmholtz}) is solved using a Krylov
subspace solver based on a ``post-conditioned'' variant of the
BiCGstab algorithm~\cite{BiCGStab} and so is replaced by
\begin{equation}
   D^{-1}_gACC^{-1}{\mathbf x}  = D^{-1}_g{\mathbf b} \hspace{0.1in}
   \equiv  \hspace{0.1in} \hat{A} CC^{-1}{\mathbf x}  = \hat{\mathbf b},
\end{equation}
where $D_g$ is the diagonal of $A$ and $C$ is the
pre-conditioner. BiCGstab is used rather than the Conjugate
Gradient algorithm~\cite{cg} due to the lack of symmetry in the derived matrix
$A$. The
application of the diagonal is performed outside the solver and
reduces the need to access this term during the application of $A{\mathbf x}$
within the Krylov method. The application of $C$ is performed through a few
(typically three) iterations of a second stationary method derived from
the Successive-Over-Relaxation method (SOR), but with a decomposition
that better reflects the nature of the atmospheric problem. This is
achieved by writing
\begin{equation}
   \hat{A} = L + T + U,
\end{equation}
where $T$ is a tridiagonal matrix arising from the vertical
discretisation, $L$ is a lower triangular matrix and $U$ is an upper
triangular matrix. Following the standard derivation of SOR,
the fixed point method is obtained.
\begin{equation}
(T + \omega L){\mathbf x}_n = (1 - \omega)T {\mathbf x}_{n-1} + \omega(r - U \mathbf{x}_{n-1} ),
\end{equation}
where $\omega = 1.5$ is the over-relaxation parameter and the only
difference from the standard method is the appearance of the
tridiagonal matrix $T$ in place of the diagonal. Note that the
lower/upper-factorisation of the matrix $T$ can be pre-computed (at each
horizontal grid point) to aid the inversion process. Furthermore,
there is no parallel decomposition in the vertical, which simplifies the
process of applying $T^{-1}$ considerably.  

\section{The Mixed-Precision Solver}
\label{sec:MPSolver}

The accuracy of the solution reached by the BiCGstab solver is
determined by the halting criterion
\begin{equation}
  \| {\bf r}_i \| = \| {\bf b} - A\cdot{\bf x}_i \|,
\end{equation}
where ${\bf x}_i$ is the solution after $i$ iterations and ${\bf r}_i$
is known as the residual vector. If the norm of ${\bf r}_i$ is small, then
${\bf x}_i$ is likely to be close to the solution ${\bf x}$. The
halting criterion for this algorithm is to break out of the iteration
loop when the residual (either absolute or relative) has fallen below
some threshold. The relative norm of the residual (R), is defined as
\begin{equation}
  \mathrm{R} = \frac{\|\mathbf{r}_i\|}{\|\mathbf{b}\|}
\label{eqn:residual}
\end{equation}
and then the halting criterion is defined as when ${\mathrm R}<{\mathrm R}_c$, where ${\mathrm R}_c$ is known
as the critical residual or the ``solver tolerance''. The optimum
value of ${\mathrm R}_c$ is
determined by a compromise between the desired {\em accuracy} of the
solution and its computational {\em cost}. 

In floating point arithmetic, the Unit of 
Least Precision (ULP) is defined in~\cite{harrison-hol99} as 
follows. The ULP in $x$ is the distance between the closest two 
{\em straddling} IEEE floating point numbers $a$ and $b$, {\em i.e.}
those with $a \le x \le b$. For numbers 
${\mathcal O}(1)$ this is ${\mathcal O}(10^{-7})$ for 32-bit numbers and 
${\mathcal O}(10^{-15})$ for 64-bit. If the stopping criterion 
\begin{equation}
\label{eq:SPlimit}
{\mathrm R}_c \gg \epsilon_{32},
\end{equation}
where $\epsilon_{32}$ is the ULP 
or machine precision for 32-bit floating point numbers, then 32-bit floating point 
numbers can satisfy the desired accuracy. Iterative methods in general are 
discussed in~\cite{linearSys}. In chapter 4, 
section 2, on the stopping criterion, the authors discuss accuracy and 
error. They show that the stopping criterion should not be set such 
that ${\mathrm R}_c < \epsilon $, Where $\epsilon$ is the precision in
the calculation.
used Conversely, decreasing $\epsilon$ from, for instance, 
$10^{-7}$ for 32-bit precision to $10^{-15}$ for 64-bit precision has no 
effect on the accuracy of the solution if ${\mathrm R} \gg \epsilon_{\rm 32}$. 

In the case of the semi-implicit time stepping scheme, the
solution to the linear system is part of the solution
procedure~\cite{QJ:QJ2235} to a larger, non-linear system and the
accuracy of the solve is dictated by the need for stability of the
time stepping scheme. Moreover, it is only an indirect
measure of the error. It is worth noting that since the finite
difference approximations to the pressure gradient are at best second
order, there is a limit to the effect of tightening the solver
tolerance on the pressure. {\em i.e.} once the error in the solver is
sufficiently small, the discretisation error becomes dominant.

In its first implementation, all operational systems using ENDGame
used a tolerance of ${\mathrm R}_c=10^{-3}$, which easily
satisfies the condition laid out in equation~(\ref{eq:SPlimit}).
These calculations can be 
done in 32-bit without loss of accuracy compared to a calculation done
in 64-bit, which offers a significant optimisation opportunity,
especially as the solver and the pre-conditioner are memory bound. A
32-bit version of the code has in effect, twice the cache memory
available compared to a 64-bit version.
This motivated the use of the mixed-precision solver, which
contributed to the optimisations that made the implementation
computationally affordable.

In principle, the whole BiCGStab algorithm could be written in
32-bit. However, the rest of the model remains 64-bit. The two
components have been combined and it is natural to encapsulate the
solver algorithm as a subroutine. It is not possible in Fortran to
coerce a 64-bit real to a 32-bit real as an argument to a subroutine,
modify its value and promote it back safely to a 64-bit real when the
subroutine exits. So, rather than pay the cost of a whole domain memory copy
from 64- to 32-bit when entering the routine and again on exit, the
pressure field is kept as a 64-bit data-type and so mixed precision
arithmetic is required. Most of the operations, including
communications, can be performed in single precision and the Fortran
run-time can coerce or promote the remaining operations automatically.

Following the initial implementation of ENDGame, numerical noise was
observed in the horizontal wind field near the pole, which could be
alleviated by tightening the solver tolerance~\citep{Walters:2017},
although it remains unclear whether the noise was a direct feature of
insufficient convergence in the vicinity of the pole, or was created by
some other numerical issue and removed/reduced through additional
iterations of the pressure solver.
The latest operational global configuration uses a tolerance of
${\mathrm R}_c=10^{-4}$, which although an order of magnitude tighter than
originally used, still satisfies the inequality in~(\ref{eq:SPlimit}).

\section{Numerical and computational performance}
\label{sec:res}
\begin{figure}[ht]
\centering\includegraphics[width=1.0\linewidth]{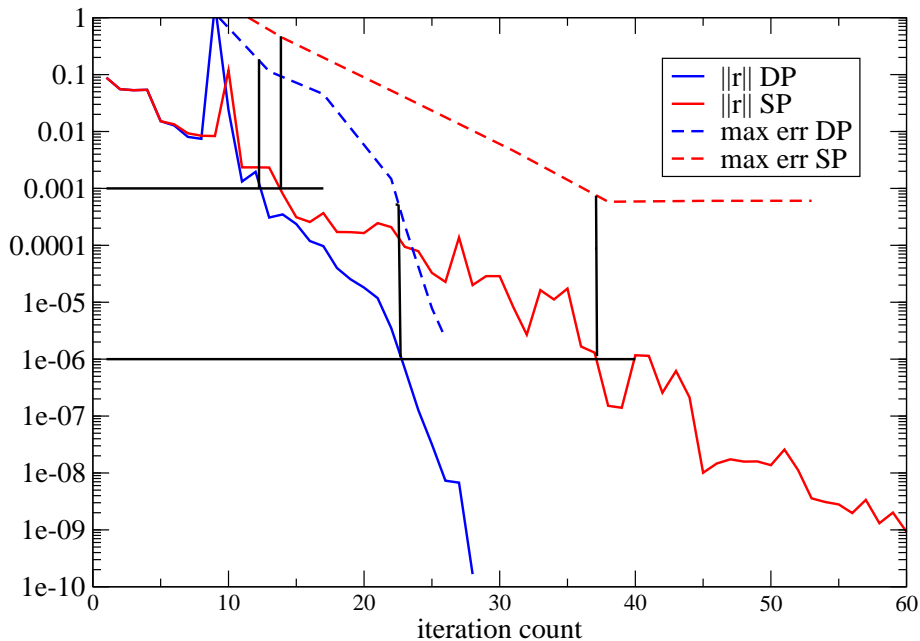}
\caption{\label{fig:residual_toy} The relative norm of residual, R 
  defined in equation~(\ref{eqn:residual}) versus iteration 
  number. The blue line, labelled ``DP'' shows data from the 64-bit
  solver, the red line, labelled ``SP'' is the 32-bit solver. 
}
\end{figure}

To test the convergence of the solver, a serial toy-model of the ENDGame
solver was constructed. The toy-model is simply the BiCGstab solver
and pre-conditioner with $\mathbf{x}$ initialised to a representative pressure
field, the Helmholtz matrix with the correct structure for
ENDGame and representative right-hand side $\mathbf{b}$ values.  No
halting criterion is used, but instead the solver is run to very high
levels of convergence, to study the convergence behaviour of R.
Shown in figure~(\ref{fig:residual_toy}) is the value of R for a 32-bit
and 64-bit solver as a function of time step. 
Both 32- and 64-bit solvers show the non-monotonic convergence typical
of the BiCGStab algorithm. For the first few iterations, convergence
is sufficiently similar to be indistinguishable. Around ${\mathrm R}=10^{-3}$,
the 32-bit version takes one extra iteration to reach the same value
of R. This ``iteration gap'' grows as the residual falls. At tighter
convergence, the 32-bit version takes twice as many iterations to
reach ${\mathrm R}=10^{-7}$ as the 64-bit version. This may well be due to the
loss of orthogonality in the Krylov sub-space vectors. In the BiCGstab
algorithm, these vectors are updated each iteration, they are not
re-computed. Accumulated round-off errors set in earlier for the lower
precision versions. The estimate of these vectors becomes worse with
increasing iteration number resulting in a loss of efficiency. This
can be recovered by restarting the algorithm; this comes at the cost of
re-computing the vectors, but this restores the orthogonality of the
Krylov sub-space. Intriguingly, the break-even point for the 32-bit
solver is ${\mathrm R}\approx 10^{-7}$, {\em i.e.} the limit of 32-bit precision. Here, the 32-bit
version takes twice as many iterations as the 64-bit version but each
iteration is twice as fast.

\begin{figure}
\centering\includegraphics[width=1.0\linewidth]{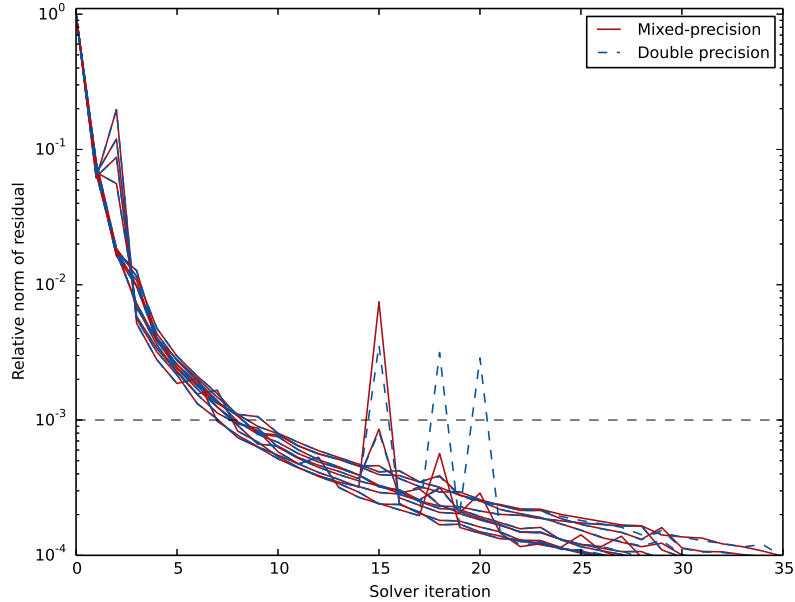}
\caption{\label{fig:residual_zoom} The relative norm of residual, R 
  defined in equation~(\ref{eqn:residual}) versus iteration 
  number from the first call to the solver in 11 simulation dates each run with
  the mixed- and double precision solver. 
}
\end{figure}

The decision to implement the mixed-precision solver in the full model
was based on tests using comprehensive weather and climate model
simulations. These included high resolution 5-day NWP forecast runs
and further testing in an operational-like research NWP
system. Results with the mixed-precision and double precision solvers
showed no noticeable scientific differences and hence the
mixed-precision solver was accepted.  To study the behaviour of the
solver in more detail in an operational-like forecast, a set of N1280 global model simulations
were run for a series of 11 dates each separated by one week, each on 
96 nodes of the Met Office Cray XC40\footnote{To improve the robustness of the timing information, each forecast was submitted 3 times, to improve the sampling of run-to-run variability. For a given date/experiment, the results of these resubmissions are identical to the bit level, so these do not add to the sampling of iteration counts or solver convergence.}. 
Mixed- and double precision solver runs for each date were 
initialised from 64-bit model states ${\mathbf x}$, saved from a previous short forecast\footnote{Operational forecast runs usually start from 32-bit truncated states for efficiency, and will usually have some special treatment ahead of and during the first time step, such as the inclusion of data assimilation increments and the use of a fully implicit first time step to better handle any imbalance these may cause. Using 64-bit dumps and switching these additional options off makes the first time step of these runs representative of later ``typical'' model time steps.}. 
Each XC40 node is dual socket, each
with an 18-core Intel Xeon (Broadwell) processor. The model is
parallelised with both MPI and OpenMP by domain decomposition across
the two dimensional horizontal domain. It was run with six MPI ranks
per socket and three OpenMP threads per MPI rank. Again, the
convergence was measured using the residual R defined in
equation~(\ref{eqn:residual}). Shown in figure~(\ref{fig:residual_zoom}) is the
value of R versus the iteration number, for the mixed- and double
precision version, up to the operational convergence criterion of
${\mathrm R}_c=10^{-4}$,
The data for convergence of the solvers is also
shown in table~(\ref{tab:timing}).

The results presented from each date are from the first call to the
solver in the first model time step, so that the prognostic fields on
entering the solver are identical between the mixed- and double
precision versions, so as to make it simple to compare the behaviour
between the two. There are several features in the figure worth remarking 
on. Firstly, for the first five iterations, the fall in the 
size of the residual is quite steep, showing quick convergence,
although four of the eleven sets of simulations demonstrate the lack
of guaranteed monotonicity, even this early on.
Secondly, up to and slightly beyond the original operational halting criterion of
${\mathrm R}_c=10^{-3}$, the values of ${\mathrm R}_{\rm MP}$ and ${\mathrm R}_{\rm DP}$ 
are very close, where $MP$ denotes mixed-precison and $DP$ double or
64-bit precsion. After a single iteration, these agree for every date to seven
significant figures, {\em i.e.} roughly the precision of a
32-bit number. Also, for each date, the residual reaches the halting
criterion in an identical number of iterations, and the values of
${\mathrm R}_{\rm MP}$ and ${\mathrm R}_{\rm DP}$ at that point all
agree to at least four significant figures.
Between ${\mathrm R}=10^{-3}$ and  ${\mathrm R}=10^{-4}$, the rate
of convergence slows down and in a few cases jumps to a value an order
of magnitude higher before dropping back down again, although these
jumps are present in both the mixed- and double precision solver, and
appear not to be related to the precision of the calculations.

The time taken to reach ${\mathrm R}\le 10^{-3}$ and ${\mathrm R}\le
10^{-4}$ is shown in Table~\ref{tab:timing}. The mixed-precision
BiCGStab is almost twice as fast per iteration than the 64-bit
version, for the reasons outlined in section~\ref{sec:intro},
{\em viz.} that the smaller data type halves the data movement both from memory to CPU and remote
communication and doubles the number of data items in the memory
caches. 
For example, the local volume for each MPI rank/OpenMP
thread is
\begin{equation}
  LV = (l_x+2h_x) \times (l_y+2h_y) \times n_{levels} / n_{threads}
\end{equation}
where $l$ is local domain size, $h$ the halo size, $x$ the East-West
direction, $y$ the North-South directions, $n_{levels}$ is the
number of vertical levels and $n_{threads}$ is the the number of
OpenMP threads the data are shared between. 
For a 96 node (3456 core) job, a typical 
decomposition might be 72 MPI ranks East-West and 16
North-South with 3 OpenMP threads resulting in a local volume of 
\begin{equation}
LV = \left(\frac{2560}{72} + 2\right) \times 
   \left(\frac{1920}{16} + 2\right) \times 
   70 / 3 = 111720
\end{equation} 
The pre-conditioner routine,
called \verb+tri_sor+, requires ten domain valued arrays for the
7-point stencil and back substitution. Thus the size of the working
set, $N$, is 
\begin{equation}
N=111720\times 10 = 1117200
\end{equation}
The memory footprint for a working set of this size is 4468800 bytes
in 32-bit and 8937600 bytes in 64-bit. The Met Office supercomputer
processors are Intel Xeon E5-2695 v4, which has a level 2 (L2) cache of
size 4718592 bytes. The working set would fit into L2 cache in
32-bit. Of course, what data are resident in cache is much more
complicated than simple size. However, the working set clearly cannot
fit into L2 cache in 64-bit. Other levels of cache, main memory
bandwidth and remote communication bandwidth will all affect the speed
of computation and by using smaller data these can be exploited more effectively.

\begin{table}
\centering
\caption{\label{tab:timing} The mean (standard error) iteration count and wall clock time taken in 
  seconds for the mixed- and double-precision solvers to reach the
  converge to the halting criterion ${\mathrm R}_c$.}
\begin{tabular}{crrrr}
  &\multicolumn{2}{c}{Mixed-precision}& \multicolumn{2}{c}{Double precision} \\
  R   & Iterations & Time (s) & Iterations & Time (s) \\\hline\hline
$10^{-3}$ & $8.4(2)$ & $0.346(6)$ & $8.4(2)$ & $0.592(9)$ \\
$10^{-4}$ & $29.3(9)$ & $1.18(2)$ & $29.3(9)$ & $2.06(3)$ \\
\end{tabular}
\end{table}

\begin{figure}
\centering\includegraphics[width=1.0\linewidth]{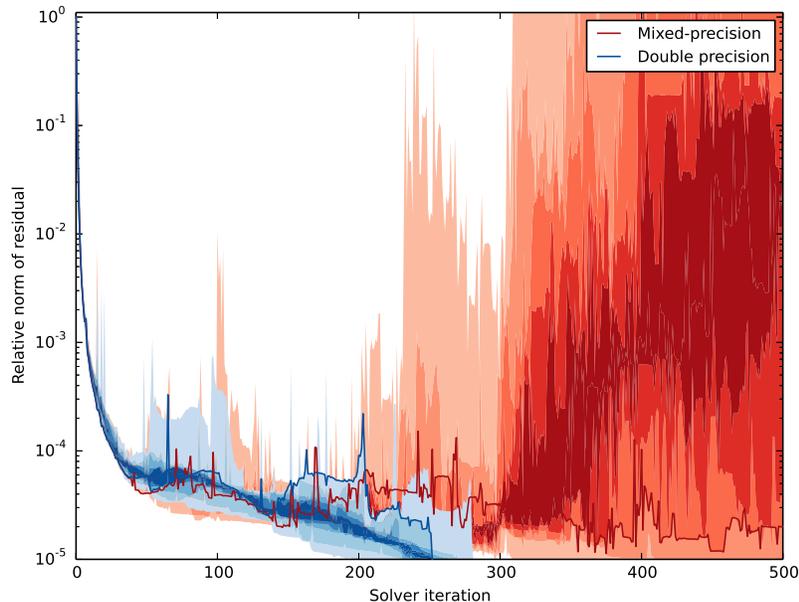}
\caption{\label{fig:residual}. The relative norm of residual, R 
  defined in equation~(\ref{eqn:residual}) versus iteration 
  number from the same simulations presented in
  figure~\ref{fig:residual_zoom}, but with the tighter halting
  criterion of ${\mathrm R}_c=10^-5$. Shading is used to highlight the
  ranges of ${\mathrm R}$ per iteration, with a ``typical'' single case
  highlighted in the solid lines.
} 
\end{figure}

To see the effect on the iteration count and convergence of the
algorithm for this work, the simulations from
figure~(\ref{fig:residual_zoom}) have been extended to use a tighter
halting criterion of ${\mathrm R}_c=10^{-5}$, which is shown in figure~(\ref{fig:residual}). 
This shows that when ${\mathrm R}<10^{-4}$, the BiCGstab algorithm starts to break
down as beyond this point, the value of the residual falls slowly and
often jumps to a higher value between iterations. The simulations with
the double precision solver take hundreds of iterations to converge to
${\mathrm R}\le10^{-5}$, and can experience many jumps before
this occurs. Below ${\mathrm R}\approx10^{-4}$, however, it is clear
that the mixed-precision approach suffers from severe errors and
after 500 iterations, only two of the simulations have converged,
whilst others have become unstable and the residual has started
growing rather than reducing with time.

In operational simulations, occasional problems with convergence have
been observed even with ${\mathrm R}_c=10^{-4}$. In these cases, the scalar
weights of the BiCGStab algorithm are close to or equal to zero and dividing by
these very small numbers can result in floating-point overflows and
the algorithm failing. The scalars are calculated by a global sum,
which in the original mixed-precision solver were performed in 32-bit,
with the local summation performed first, then an MPI reduction call
on a single number. To protect against the scalars summing to
zero, however, the global sum only has since been reverted to 64-bit. This has a
negligible cost as the MPI reduction is the dominant cost of the
procedure and this is latency bound. Reverting the global sums to
64-bit reduces the likelihood of the model failing with overflowing
fields, but the problem of slow convergence remains.  Slow convergence
can be the result of a loss of orthogonality of the Krylov sub-space,
although there are multiple reasons for which BiCGStab can break
down~\cite{Graves-Morris2002}. To address this and improve the
robustness of operational systems, the implementation of BiCGStab was
modified to include a mandatory restart if convergence has not been
achieved after a fixed number of iterations.
The restart threshold is set at a relatively high number of
iterations of 150.
The results above suggest that this poor convergence could be either due
to the general convergence issues seen in both the mixed-
and double precision solvers, or due to the additional rounding errors
seen in the mixed-precision solver with very small values of $R$.

\section{Conclusion}
\label{sec:con}

Speed, accuracy and stability of computation are all important criteria for numerical
calculations, particularly for an operational NWP system. When
carefully controlled, the use of reduced precision  can have a
big impact on the speed of a calculation without affecting its
accuracy or stability. However, as
shown in this work, for a complex system that is susceptible to
non-convergence of the numerical algorithm employed to solve it,
reduced precision can only be safely employed within a certain regime.

There are undoubted benefits to reducing precision, in this case
halving the run-time of the solver. However, other numerical
instabilities such as the polar noise described in a previous section
have forced the model to be run on a regime where the iterative,
Krylov solver has problems with slow convergence. This makes the use
of reduced precision all the more important as the run-times with a
tighter stopping criteria are longer. Naively, it should be possible
to run with reduced precision at tighter stopping condition than
$10^{-4}$, however, problems with slow convergence and numerical
stability indicate that the algorithm used, BiCGStab, is failing,
either to converge or doing so far too slowly. This is true for the
double-precision version and so it is safe to conclude that these
issues are not caused by using reduced precision. However, as the
reduced precision version appears to suffer more severe problems,
including failing to converge, it is clear that there is a complex
interplay between accumulated round-off errors and other errors
introduced by the numerical scheme.

These issues, combined with the behaviour shown in
figure~(\ref{fig:residual}), suggest that there could be some benefit
from continuing to investigate the cause of the noise in the
horizontal wind field near the pole. This may allow the operational forecasts to revert
to the slacker solver tolerance of $R_c=1.0\times 10^{-3}$ and hence
make the solution less susceptible to problems with convergence,
whilst reducing computational cost.
These results also make clear that the poor convergence is not directly 
related to the mixed-precision solver, and would still be
present if the solver were routinely run solely in double
precision. This should encourage atmospheric model developers to
continue to consider reducing the precision of some calculations to
improve the efficiency of their simulations.

\section{acknowledgements}
   The authors thank T. Allen and N. Wood, both from the Met Office,
   for in-depth discussions and advice.


\end{document}